\documentclass[12pt,a4paper,titlepage]{article}

\usepackage{amsfonts}
\usepackage{amsmath}
\usepackage{amssymb}
\usepackage{amsthm}
\usepackage{graphicx} 
\usepackage{float}
\usepackage{caption}
\usepackage{tikz-cd}
\usepackage[utf8]{inputenc}
\usepackage{multicol}
\usepackage{xcolor}

\newtheorem*{twr*}{Theorem}
\newtheorem{twr}{Theorem}[section]
\newtheorem{wn}[twr]{Corollary}
\newtheorem{uwa}[twr]{Remark}
\newtheorem{lem}[twr]{Lemma}
\newtheorem{defi}[twr]{Definition}
\newtheorem{przy}[twr]{Example}

\newcommand{\doesnotdivide}{\hspace{-5pt}\not\vert}

\usepackage[T1]{fontenc}
\usepackage[top=2cm, bottom=2cm, left=3cm, right=3cm]{geometry}
\makeatletter
\title{\huge The jump of the Milnor number of quasihomogeneous singularities for linear deformations}
\author{{\huge Aleksandra Zakrzewska}\\Faculty of Mathematics and Computer Science\\ University of Lodz\\ul. Banacha 22\\90-238 Lodz, Poland\\aleksandra.zakrzewska@wmii.uni.lodz.pl\\
	ORCID:0000-0002-7523-5133
}
\date{}
\begin{document}
%\begin{center}
%{\Large \textsc{\textbf{The jump of Milnor number of quasihomogeneous singularities for linear deformations}}}
%\vspace{1cm}
%\\
%\textsc{Aleksandra Zakrzewska}
%\vspace{1cm}
%\\
%\end{center} 
\maketitle
\begin{abstract}
The jump of the Milnor number of an isolated singularity $f_0$ is the minimal non-zero difference between the Milnor numbers of $f_0$ and one of its deformations
$f_s$. We determinate the jump of quasihomogeneous singularities
in the class of linear deformations.
\end{abstract}

\section{Introduction}
One of the important problems in singularity theory is the adjacency problem: when a singularity (or a class of singularities) can be deformed to another one. In other words whether a "type" of a singularity may be changed to another "type" be an arbitrarily small deformation. A simpler problem is to find how some invariants of singularities may change by an arbitrarily small deformation. In the article we study such a change of the Milnor number for isolated plane curve singularities. We are interested in finding the smallest positive change under some class of deformations -- we will call it the jump of the Milnor number of a given singularity.

We start from basic definitions. They are given in $n$-dimensional case, but further we will focus on only the \textbf{plane curve singularities}. Let $f_0:\left(\mathbb{C}^n,0\right)\to\left(\mathbb{C},0\right)$ be an  \textbf{isolated singularity} or in short \textbf{singularity}. We define a \textbf{deformation of the singularity} $f_0$ as a germ of a holomorphic function  $f:\left(\mathbb{C}\times\mathbb{C}^n,0\right)\to\left(\mathbb{C},0\right)$ such that
\begin{enumerate}
	\item $f(0,z)=f_0(z)$,
	\item $f(s,0)=0$.
\end{enumerate}

The deformation $f(s,z)$ of the singularity $f_0$ will be treated as a family $(f_s)$ of function germs, taking $f_s(z):=f(s,z)$. For the sufficiently small $s$ we can define the \textbf{Milnor number of $f_s$ at $0$} by $$\mu_s:=\mu(f_s)=\dim_{\mathbb{C}}{\mathcal{O}_n}/{\left(\nabla f_s\right)},$$
where $\mathcal{O}_n$ is the ring of holomorphic function germs at $0$, and $(\nabla f_s)$ is the ideal in $\mathcal{O}_n$ generated by $\frac{\partial f_s}{\partial z_1},\ldots,\frac{\partial f_s}{\partial z_n}$.

The Milnor number is upper semi-continuous in the Zariski topology in families of
singularities (\cite{GLS}, Theorem 2.6 I and Proposition 2.57 II), so there exists an open neighbourhood $0\in S$ such that
\begin{enumerate}
	\item $\mu_s=\textrm{const.}$ for $s\in S\setminus \{0\}$,
	\item $\mu_0\geq \mu_s$ for $s\in S$.
\end{enumerate}

The constant difference $\mu_0-\mu_s$ (for $s\in S$) will be called \textbf{the jump of the deformation $(f_s)$} and denoted by $\lambda((f_s))$. \textbf{The jump of the Milnor number of the singularity} $f_0$ is the smallest non-zero value among all the jumps of deformations of the singularity $f_0$. It will be denoted by $\lambda(f_0)$.

Many authors have considered what values the jump of the Milnor number can take. One of the first general result was obtained by Sabir Gusein-Zade (\cite{SGZ}). In his work he proved that there exist singularities $f_0$ for which $\lambda(f_0)> 1$ and that for any irreducible plane curve singularity $f_0$ we have $\lambda(f_0)=1$. Later, S. Brzostowski, T. Krasi\'nski and J. Walewska in \cite{BKW} proved that for the particular reducible singularities $f^n_0(x,y)=x^n+y^n$, $n\geq 2$, we have $\lambda(f_0)=\left[\frac{n}{2}\right]$. Determining the jump of a singularity is difficult because it is not a topological invariant (\cite{BK}, \cite{PW} Section 7.3).

A simpler problem is to determinate the jump when we limit ourselves to specific classes of deformations. For non-degenerate deformations (it means each element of the family $(f_s)$ is a non-degenerate singularity in the Kouchnirenko sense \cite{KU}) the jump (denoted by $\lambda^{nd}(f_0)$)  was considered in \cite{BOD}, \cite{WAL}, \cite{BKW}, \cite{KW}.

In this paper we consider the jump of the Milnor number for \textbf{linear deformations} of $f_0$ i.e. deformations of the form $f_s=f_0+sg$, where $g$ is  a holomorphic function in the neighbourhood of $0$ such that $g(0)=0$. We will denote the jump of $f_0$ for this class of deformations by $\lambda^{lin}(f_0)$. The main result is a formula for the jump of the Milnor number $\lambda^{lin}(f_0)$ for quasihomogeneous plane curve singularities. The simpler problem of homogeneous singularities was treated in \cite{ZAK}.

In generic case (the general precise result is given in Theorem \ref{quasitwr}) the formula is as follows
\begin{twr*}
	If $f_0(x,y)=a_{p,0}x^p+\ldots +a_{0,q}y^q$ is a quasihomogeneous isolated singularity with generic coefficients and $3\leq p\leq q$ then
	$$\lambda^{lin}(f_0)=\left\{\begin{array}{ll}
		p-2,&\textrm{if }p=q\\
		p-1,&\textrm{if }p\neq q\textrm{ and }p\vert q\\
		GCD(p,q),&\textrm{if }p\neq q\textrm{ and }p\doesnotdivide q
		\end{array}\right..$$
\end{twr*}
The first case concerns the homogeneous singularity. We illustrate the result with two examples.
\begin{przy}
	For a homogeneous singularities $f_0(x,y)=x^n+y^n$, where $n\geq 3$, the various types of jumps are different:
	$$\lambda(f_0)=\left[\frac{n}{2}\right],
	\lambda^{lin}(f_0)=n-2,
	\lambda^{nd}(f_0)=n-1.
	$$
	If we put for example $n=5$ then:
	$$\lambda(f_0)=2,\lambda^{lin}(f_0)=3,\lambda^{nd}(f_0)=4.$$
\end{przy}  
\begin{przy}\label{prz2}
	For the quasihomogeneous singularity $f_0(x,y)=x^6+y^9$ we have $\lambda^{nd}(f_0)=\lambda^{lin}(f_0)=3$ but the constructions given in \cite{WAL} for non-degenerate case, and in Theorem \ref{latw} for linear case give different deformations realizing this jump:
	\begin{enumerate}
		\item $f_s(x,y)=x^6+y^9+sx^5y$ -- the non-degenerate deformation,
		\item $f_s(x,y)=x^6+y^9+sxy(y^3+x^2)^2$ -- the linear deformation.
	\end{enumerate}
\end{przy}
To get the main the \textbf{Enriques diagrams} will be used. To any singularity we  assign a weighted Enriques diagram $(D,\nu)$ which represents the whole resolution process of this singularity (\cite{ECA} Chapter 3.9). It is a tree with two types of edges. M. Alberich-Carrami\~{n}ana and J. Ro\'e (\cite{MAC} Theorem 1.3, Remark 1.4) gave a necessary and sufficient condition for two Enriques diagrams of singularities to be linear adjacent. It means that one singularity is a linear deformation of another. They used a wider class of Enriques diagrams, so-called abstract Enriques diagrams, which are described in Section 2.

\section{Abstract Enriques diagrams}
Information about abstract Enriques diagrams can be found in \cite{MAC} and \cite{KP}. Moreover in my previous paper \cite{ZAK}, in which I gave the estimation of $\lambda^{lin}(f_0)$ for homogeneous singularities, abstract Enriques diagrams are described in more details with examples. The formula for $\lambda^{lin}$ for homogeneous singularities is in my PhD thesis \cite{doktorat} (in Polish).

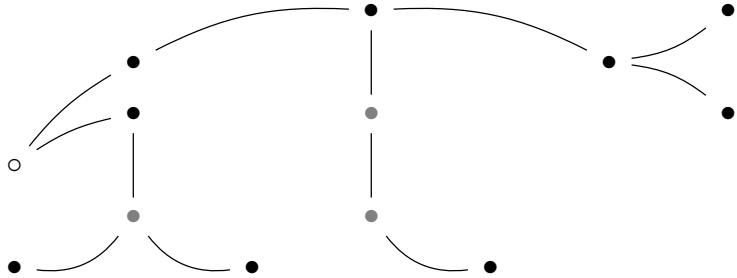
\begin{figure}[ht]
	\begin{center}
		\begin{tikzcd}[row sep=tiny]
			&&&\bullet\arrow[drr,dash,bend left=15]\arrow[dd,dash]&&&\bullet\\
			&\bullet\arrow[urr,dash,bend left=15]&&&&\bullet\arrow[dr,dash,bend left=15]\arrow[ur,dash,bend right=15]&\\
			&\bullet\arrow[dd,dash]&&{\color{gray}\bullet}\arrow[dd,dash]&&&\bullet\\
			\circ\arrow[uur,dash,bend left=10]\arrow[ur,dash,bend left=10]&&&&&&\\
			&{\color{gray}\bullet}\arrow[dr,dash,bend right]\arrow[dl,dash,bend left]&&{\color{gray}\bullet}\arrow[dr,dash, bend right]&&&\\
			\bullet&&\bullet&&\bullet&&\\	
		\end{tikzcd}
	\end{center}
	\caption{\label{rys1B}The abstract Enriques diagram. Satellite vertices are marked in gray. The root is white.}
\end{figure}

\begin{defi}[\cite{MAC}]
	An \textbf{abstract Enriques diagram} (in short an Enriques diagram) is a rooted tree $D$ with binary relation between vertices, called proximity, which satisfies:
	\begin{enumerate}
		\item The root is proximate to no vertex.
		\item Every vertex that is not the root is proximate to its immediate predecessor.
		\item No vertex is proximate to more than two vertices.
		\item If a vertex $Q$ is proximate to two vertices, then one of them is the immediate predecessor of $Q$ and it is proximate to the other.
		\item Given two vertices $P,Q$ with $Q$ proximate to $P$, there is at most one vertex proximate to both of them.
	\end{enumerate}
\end{defi}
The fact that $Q$ is proximate to $P$ we will denote by $Q\to P$. The vertices which are proximate to two points are called \textbf{satellite}, the other vertices (except the root) are called \textbf{free}. The vertex is \textbf{final} if it has no successor. To show graphically the proximity relation, Enriques diagrams are drawn according to the following rules:
\begin{enumerate}
	\item If $Q$ is a free successor of $P$, then the edge going from $P$ to $Q$ is curved.
	\item The sequence of edges connecting a maximal succession of vertices proximate to the same vertex $P$ are shaped into a line segment, orthogonal to the edge joining $P$ to the first vertex of the sequence (if it is also straight).
\end{enumerate}
The example of an abstract Enriques diagram is shown in Figure \ref{rys1B}.

We will now introduce few basic notations that are needed in later chapters. First, we define weights on vertices of an abstract Enriques diagrams which correspond, in particular case of plane curve singularities, to the orders of the proper transforms of the function describing the singularity.

A \textbf{weight function} is any function $\nu:D\to\mathbb{Z}$. A pair $(D, \nu)$, where $D$ is an abstract Enriques diagram and $\nu$ a weight function, is called a \textbf{weighted Enriques diagram}. A \textbf{consistent Enriques diagram} is a weighted Enriques diagram such that for all $P\in D$ 
\begin{equation}\label{cons}
	\nu(P)\geq\sum_{Q\to P}\nu(Q).
\end{equation}
A \textbf{complete Enriques diagram} is a weighted Enriques diagram such that for all non-final $P\in D$ the equality in (\ref{cons}) holds and for all final $P\in D$ it is a free vertex with weight $1$ not proximate to another free vertex with weight $1$. To the weight function $\nu$ of a weighted diagram $D$ we associate a \textbf{system of values} on $D$, which is another map ${\normalfont\textrm{ord}}_\nu:D\to\mathbb{Z}$, defined recursively as
$${\normalfont\textrm{ord}}_\nu(P):=\left\{\begin{array}{ll}\nu(P),&\textrm{if }P\textrm{ is the root,}\\\nu(P)+\sum\limits_{P\to Q}{\normalfont\textrm{ord}}_\nu(Q),&\textrm{otherwise.}\end{array}
\right.$$ For any consistent $(D,\nu)$ we define \textbf{the Milnor number of $(D,\nu)$} by
$$\mu((D,\nu)):=\sum_{P\in D}\nu(P)(\nu(P)-1)+1-r_D,$$
where $r_D:=\sum_{P\in D}r_D(P)$, $r_D(P):=\left(\nu(P)-\sum_{Q\to P}\nu(Q)\right)$ for every $P\in D$.

A \textbf{subdiagram} of an abstract Enriques diagram $D$ is a subtree $D_0\subset D$ with the same proximity relation such that if $Q\in D_0$ then its predecessor belongs to $D_0$.

In the class of weighted Enriques diagrams, we introduce equivalence relation. We say that weighted diagrams $(D,\nu)$ and $(D',\nu')$ are \textbf{equivalent} if they differ at most in free vertices of weight $1$. The equivalence class of $(D,\nu)$ is denoted by $[(D,\nu)]$ and called the \textbf{type} of $(D,\nu)$. Of course, the Milnor number is invariant in the class $[(D,\nu)]$.

A \textbf{minimal Enriques diagram} is a consistent Enriques diagram $(D,\nu)$ with:
\begin{enumerate}
	\item no free vertices of weight $0$,
	\item no free vertices of weight $1$ except for these such $P\in D$ for which there exists a satellite vertex $Q\in D$ satisfying $Q\to P$.
\end{enumerate}

It is easy to see (\cite{ZAK}, Theorem 2.12) that
\begin{twr}
	Let $(D,\nu)$ be a consistent weighted diagram. There exists exactly one minimal diagram which belongs to $[(D,\nu)]$.
\end{twr}

The theory of Enriques diagrams has its roots in the theory of plane curve singularities. The embedded resolution of a plane curve singularity using blow-ups can be explicitly presented as a complete Enriques diagram. A precise description can be found in \cite{ECA} Chapter 3.8 and Chapter 3.9. Two plane curve singularities are topologically equivalent if and only if their Enriques diagrams are isomorphic (as graphs). For the Enriques diagram of a plane curve singularity, the weight function represents the orders of the consecutive proper transforms while the system of values -- the orders of the total transforms of the function defining the singularity. Also the Milnor number of the Enriques diagram coincides with the Milnor number of the corresponding singularity. We need only the next fact which easily follows from these results.

\begin{twr}[\cite{ECA} Theorem 3.8.6]\label{bije}
	There exists a bijection between minimal Enriques diagrams and topological types of singularities.
\end{twr}

In the paper \cite{MAC}, M. Alberich-Carrami\~{n}ana and J. Ro\'e gave a necessary and sufficient condition for two Enriques diagrams of singularities to be linear adjacent. This is the key result we will use in the sequel. First we give definitions.

\begin{defi}
	Let $(D,\nu)$ and $(D',\nu')$ be weighted Enriques diagrams, with $(D',\nu')$ consistent. We will write $(D',\nu')\geq (D,\nu)$ when there exist isomorphic subdiagrams $D_0\subset D$, $D_0'\subset D'$ with an isomorphism (that preserves proximity relations) $$i:D_0\to D_0'$$ such that the new weight function $\kappa:D\to\mathbb{Z}$ for $D$, defined by
	$$\kappa(P):=\left\{\begin{array}{cr}
		\nu'(i(P)),&P\in D_0\\0,&P\notin D_0
	\end{array}
	\right.$$
	satisfies $${\normalfont\textrm{ord}}_\nu(P)\leq{\normalfont\textrm{ord}}_\kappa(P)$$ for any $P\in D$.
\end{defi}

\begin{defi}\label{def1}
	Let $[(D,\nu)]$ and $[(\widetilde{D},\widetilde{\nu})]$ be types of Enriques diagrams. $[(\widetilde{D},\widetilde{\nu})]$ is \textbf{linear adjacent} to $[(D,\nu)]$ if there exists a consistent Enriques diagram $(D',\nu')\in[(\widetilde{D},\widetilde{\nu})]$ such that $(D',\nu')\geq (D_{\textrm{min}},\nu_{\textrm{min}})$, where $(D_{\textrm{min}},\nu_{\textrm{min}})$ is the minimal diagram of type $[(D,\nu)]$.
\end{defi}

\begin{twr}[\cite{MAC} Theorem 1.3 and Remark 1.4]\label{twrmac}
	Let $[(D,\nu)]$ and $[(\widetilde{D},\widetilde{\nu})]$ be types of consistent Enriques diagrams. The following conditions are equivalent:
	\begin{enumerate}
		\item $[(\widetilde{D},\widetilde{\nu})]$ is linear adjacent to $[(D,\nu)]$.
		\item For every singularity $f_0$ whose Enriques diagram belongs to $[(\widetilde{D},\widetilde{\nu})]$, there exists a linear deformation $(f_s)$ of $f_0$ such that the Enriques diagram of a generic element $f_s$ belongs to $[(D,\nu)]$.
		\item There exists a singularity $f_0$ whose Enriques diagram belongs to $[(\widetilde{D},\widetilde{\nu})]$ and a linear deformation $(f_s)$ of $f_0$ such that the Enriques diagram of a generic element $f_s$ belongs to $[(D,\nu)]$.		
	\end{enumerate}
\end{twr}
This theorem was also formulated using prime divisors by J. Fernández de Bobadilla, M. Pe Pereira and P. Popescu-Pampu in Theorem 3.25 (\cite{BPP}).

Theorems \ref{bije} and \ref{twrmac} imply the following corollary:
\begin{wn}
	$\lambda^{lin}(f_0)$ is a topological invariant.
\end{wn}

\section{Enriques diagrams of quasihomogeneous singularities}
Let $f_0(x,y)=\sum_{i,j\in\mathbb{N}}a_{i,j}x^iy^j$ be an isolated singularity. It is known that $f_0$ is reduced in the ring $\mathbb{C}\{x,y\}$ of convergent series. The singularity $f_0$ is called \textbf{quasihomogeneous}, if there exist $w_x,w_y\in\mathbb{N}$ and a number $W\in\mathbb{N}$ such that, for every $(i,j)\in\textrm{supp}(f_0)$, it holds $iw_x+jw_y=W$, where $\textrm{supp}(f_0):=\{(i,j)\in\mathbb{N}:a_{i,j}\neq 0\}$. Without loss of generality, $f_0$ can be expressed as
\begin{equation}
	f_0(x,y)=x^ky^l(x^p+\ldots+\gamma_{i,j}x^iy^j+\ldots+\gamma_{0,q}y^q),\quad k,l\in\{0,1\},p\leq q,k+l+p\geq 2,
\end{equation}
and for every term $\gamma_{i,j}x^iy^j$, $\gamma_{i,j}\neq 0$, the equality $(i+k)w_x+(j+l)w_y=W$ holds. 

Then after simple rescaling the variables $x\mapsto x', y\mapsto \frac{y'}{\sqrt[q+l]{\gamma_{0,q}}}$, that does not change the Milnor number of $f_0$, we may assume $f_0$ has the form:
\begin{equation}\label{quasi}
	f_0(x,y)=x^ky^l(x^p+\ldots+\gamma_{i,j}x^iy^j+\ldots+y^q),\quad k,l\in\{0,1\},p\leq q,k+l+p\geq 2,
\end{equation}
In the case $p=q$ we get a homogeneous singularity.

Since $f_0$ is reduced and quasihomogeneous in two variables, we can represent $f_0$ as a product of irreducible factors
\begin{equation}\label{product}
	f_0(x,y)=x^ky^l\prod_{i=1}^{\tilde{d}}\left(x^r+\alpha_iy^s\right),\quad\alpha_i\neq 0,\alpha_i\neq\alpha_j\textrm{ for }i\neq j,
\end{equation}
where $\tilde{d}=\textrm{GCD}(p,q)$, $r=\frac{p}{\tilde{d}}$, $s=\frac{q}{\tilde{d}}$, $\textrm{GCD}(r,s)=1$. By this form of quasihomogeneous singularity and by the resolution process of singularities (more details in \cite{ECA} Chapter 3.7) the Enriques diagram of any quasihomogeneous singularity can be easily described.

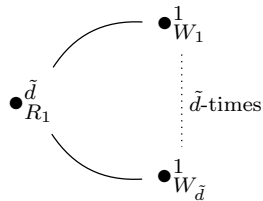
\begin{figure}[ht]
	\begin{center}
		\begin{tikzcd}[row sep=tiny]
			& \bullet^1_{W_1} \arrow[dd,dash,dotted,"\tilde{d}\textrm{-times}"] \\
			\bullet^{\tilde{d}}_{R_1} \arrow[ur, dash, bend left] \arrow[dr, dash, bend right] & \\
			& \bullet^1_{W_{\tilde{d}}}
		\end{tikzcd}
	\end{center}
	\caption{The Enriques diagram of a homogeneous singularity of order $\tilde{d}$.\label{fhomo}}
\end{figure}

\begin{figure}[ht]
		\begin{center}
			\begin{tikzcd}[row sep=tiny]
				&&&&&\bullet^1_{W_1} \arrow[dd,dash,dotted,"\tilde{d}\textrm{-times}"]\\
				\bullet_{R_1}^{\nu(R_1)}\arrow[dr, dash, bend right]&&\bullet_{R_m}^{\nu(R_m)}\arrow[dr,dash,dotted]&&\bullet^{\tilde{d}}_{R_t} \arrow[ur, dash, bend left] \arrow[dr, dash, bend right]&  \\
				&\bullet_{R_2}^{\nu(R_2)}\arrow[ru,dash,dotted,bend right]&&\bullet\arrow[ur,dash,dotted]&&\bullet^1_{W_{\tilde{d}}}
			\end{tikzcd}
		\end{center}
		\caption{The Enriques diagram of a quasihomogeneous singularity $f_0$ for $k=l=0$.\label{ftilde}}
	\end{figure}
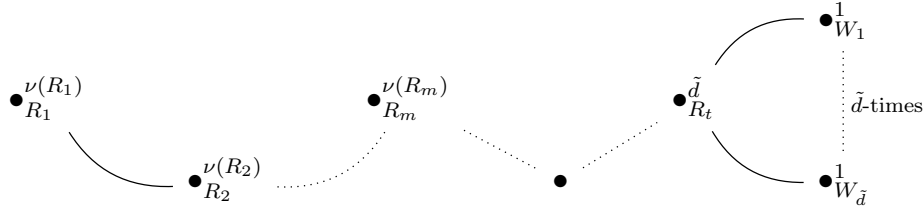

In fact, let assume first that $k=l=0$. If $r=s$ then singularity (\ref{product}) is homogeneous and hence $r=s=1$ and $p=q=\tilde{d}$. So $f_0(x,y)=\prod_{i=1}^{\tilde{d}}\left(x+\alpha_iy\right)$ for some $\alpha_i\neq 0,\alpha_i\neq\alpha_j\textrm{ for }i\neq j$. Then one blowing up resolves the singularity and the Enriques diagram of $f_0$ is shown in Figure \ref{fhomo}. Now assume $r<s$ (the case $s<r$ is analogous). So $f_0(x,y)=\prod_{i=1}^{\tilde{d}}\left(x^r+\alpha_iy^s\right)$, $r<s$, $\textrm{GCD}(r,s)=1$. Hence the singularity $f_0$ has the unique tangent line $\{x=0\}$. Then after one blowing up the proper transform of this singularity is described in the coordinates $(x',y')=(\frac{x}{y},y)$ by the polynomial $\prod_{i=1}^{\tilde{d}}\left(x'^r+\alpha_iy'^{s-r}\right)$. This singularity has also the unique tangent line (either $\{x=0\}$ if $r<s-r$ or $\{y=0\}$ if $r>s-r$) except the case $r=1$ and $s=2$. In the exceptional case we get a homogeneous singularity. In the first case (only one tangent line) after finite number of blowing ups we also get a homogenous singularity. In both cases we always get a homogenous singularity for which the next blowing up gives its resolution. According to the above description we may describe the Enriques diagram $(D,\nu)$ of $f_0$ (see Figure \ref{ftilde}). The first edges (from $R_1$ to some $R_m$) are curved and next ones (from $R_m$ to $R_t$) are straight. The diagram $(D,\nu)$ has $\tilde{d}$ final vertices. Moreover this is a complete Enriques diagram. If $p\vert q$ then $t=\frac{q}{p}$. In particular if $f_0$ is homogeneous then $t=1$.

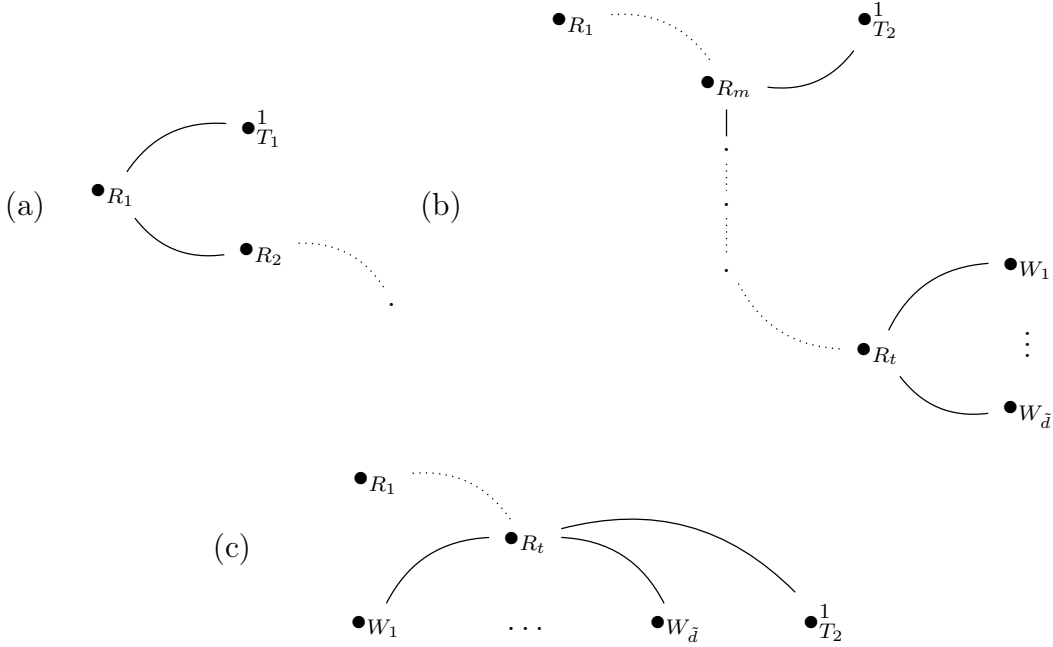
\begin{figure}[htp]
	\centering
	
	(a)\quad\begin{tikzcd}[row sep=tiny]
		&\bullet^1_{T_1}&\\
		\bullet_{R_1}\arrow[dr, dash, bend right]\arrow[ur, dash, bend left]&&\\
		&\bullet_{R_2}\arrow[rd,dash,dotted,bend left]& \\
		&&. \\
	\end{tikzcd}				
	(b)\begin{tikzcd}[row sep=tiny]
		&\bullet_{R_1}\arrow[rd,dash,dotted,bend left]&&\bullet_{T_2}^1&& \\
		&&\bullet_{R_{m}}\arrow[dd,dash]\arrow[ur, dash, bend right]&&& \\
		&&&&&\\
		&&.\arrow[dd,dash,dotted]&&&\\
		&&&&&\\
		&&.\arrow[dd,dash,dotted]&&&\\
		&&&&&\\
		&&.\arrow[dr, dash, bend right,dotted]&&\bullet_{W_1}&\\
		&&&\bullet_{R_t}\arrow[ur,dash,bend left]\arrow[dr,dash,bend right]&\vdots&\\
		&&&&\bullet_{W_{\tilde{d}}}&
	\end{tikzcd}
	\medskip
	\\(c)	\begin{tikzcd}[row sep=tiny]
		&\bullet_{R_1}\arrow[rd,dash,dotted,bend left]&&&& \\
		&&\bullet_{R_t}\arrow[ddr,dash,bend left]\arrow[ddrr,dash,bend left]\arrow[ddl,dash,bend right]&&&\\
		&&&&&\\
		&\bullet_{W_1}&\ldots&\bullet_{W_{\tilde{d}}}&\bullet_{T_2}^1&
	\end{tikzcd}	
	\caption{The Enriques diagrams of quasihomogeneous singularities. In the figure (a) $l=1$, while in (b) and (c) $k=1$. Case (c) holds if $p\vert q$. \label{typII}}
\end{figure}

If $k=1$ or $l=1$ then we proceed analogously as above with small modification. We have to add one or two leaves to the Enriques diagram in Figure \ref{ftilde} to appropriate vertices. If $l=1$ i.e. there is the factor $y$ in the factorization (\ref{product}) of $f_0$, we add a leaf $T_1$ with weight $1$ to the root $R_1$ (Figure \ref{typII}(a)). If $k=1$ i.e. there is the factor $x$ in the factorization (\ref{product}) of $f_0$, we add such a leaf $T_2$ to the last free vertex among $R_1,\ldots,R_t$ i.e. to $R_m$ in Figure \ref{ftilde}. Two possible cases $R_m\neq R_t$ and $R_m=R_t$ are presented in Figure \ref{typII}(b) and \ref{typII}(c), respectively. 

For $t,d\in\mathbb{N}$ we define the set $H_d^t$ as the set of the abstract Enriques diagrams $(D,\nu)$ satisfying conditions:
\begin{enumerate}
	\item $(D,\nu)$ is a minimal diagram,
	\item the elements of $D$ is a sequence $\{R_1,\ldots,R_t\}$  such that $R_i$ is a successor of $R_{i-1}$ for $i\in\{2,\ldots,t\}$ (a bamboo from $R_1$ to $R_t$),
	\item $\nu(R_t)=d$.
\end{enumerate} 

From the above construction of the Enriques diagrams of a quasihomogeneous singularity (\ref{product}) we see that its minimal diagram belongs to some $H_d^t$. We denote the subset of $H_d^t$ corresponding to quasihomogeneous singularities by $Q_d^t$. This means for every diagram $(D,\nu)$ from $Q_d^t$ there exists a singularity (\ref{quasi}) such that $d=GCD(p,q)=\tilde{d}$ (if $p$ does not divide $q$) and $d=GCD(p,q)+k=\tilde{d}+k$ (if $p$ divides $q$) and $(D,\nu)$ has the same type as the Enriques diagram of (\ref{quasi}). In particular, for $t=1$ the set $Q^t_d$ represents homogeneous singularities and for $d=1$ -- irreducible ones (omitting factors $x^ky^l$ in (\ref{quasi})).

It is easy to show the abstract Enriques diagrams which belong to $Q_d^t$ have the following properties.

\begin{twr}
	If a weighted Enriques diagram $(D,\nu)$ belongs to $Q_d^t$ ($t\neq 1$) then
	\begin{enumerate}
		\item $\nu(R_1)\leq\sum\limits_{R_i\to R_1}\nu(R_i)+1$,
		\item if $R_k$ is the first satellite vertex for some $k\in{2,\ldots,t}$ then\\ $\nu(R_{k-1})\leq\sum\limits_{R_i\to R_{k-1}}\nu(R_i)+1$
		\item for any $k=2,\ldots,t$ such that $R_{k+1}$ is not the first satellite vertex, we have $\nu(R_{k})=\sum\limits_{R_i\to R_{k}}\nu(R_i)$.
	\end{enumerate}
\end{twr} 

The subset $Q_d^t$ is a proper subset of $H_d^t$, for example the minimal Enriques diagram of the singularity  $f_0(x,y)=(x^2-y^2)(x^6-y^9)$ belongs to $H_d^t\setminus Q_d^t$.

For any $(D,\nu)\in H_d^t$ we define $w_D$ as the number of vertices which $R_t$ is proximate to. If $(D,\nu)$ is the Enriques diagram of singularity (\ref{product}), then obviously \begin{equation}\label{w}w_D=\left\{\begin{array}{ll}0,&\mathrm{if}\hspace{1ex}p=q\\1,&\mathrm{if}\hspace{1ex}p\neq q\textrm{ and }p\vert q\\2,&\mathrm{if}\hspace{1ex}p\neq q\textrm{ and }p\doesnotdivide q\end{array}\right..
\end{equation} 

\section{Estimation of the Milnor number for abstract Enriques diagrams}
In this section we will estimate the Milnor number of these diagrams to which  diagrams from $Q_d^t$ are linear adjacent. Precisely, for any $(D,\nu)\in Q^t_d$ we will find the maximum in the set \begin{equation}\label{set}
	\{\mu((E,\lambda)):[(D,\nu)]\textrm{ is linear adjacent to }[(E,\lambda)],(E,\lambda)\notin[(D,\nu)]\},
\end{equation} where $d,t\in\mathbb{N}$ and $dt>1$. If $dt=1$ then $H_1^1=Q_1^1$ represents a smooth curve (by our definition it is not a singularity). We will show that this maximum equals 
$$\begin{array}{ll}
	\mu((D,\nu))-1,&\mathrm{if}\hspace{1ex} d=1\\
	\mu((D,\nu))-1,&\mathrm{if}\hspace{1ex} d=2,w_D=0\\
	\mu((D,\nu))-w_D,&\mathrm{if}\hspace{1ex}d=2,w_D\neq 0\\
	\mu((D,\nu))-(d-2+w_D),&\mathrm{if}\hspace{1ex}d\geq 3\end{array}.$$

\begin{figure}[ht]
	\begin{center}
		\begin{tikzcd}[row sep=tiny]
			\bullet_{R_1}^{\nu(R_1)}\arrow[dr, dash, bend right]&&&  \\
			&\bullet_{R_2}^{\nu(R_2)}\arrow[r,dash,dotted]&\bullet^{\nu(R_{t-1})}_{R_{t-1}}\arrow[r,dash,dotted]&\bullet^d_{R_t} \\
		\end{tikzcd}
	\end{center}
	\caption{The minimal Enriques diagram $(D,\nu)$.\label{diagrD}}
\end{figure}
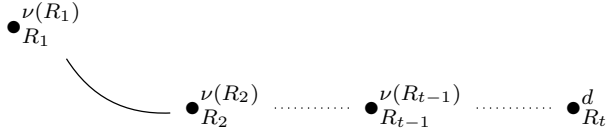

\begin{figure}[htp]
	
	\qquad(a)\begin{tikzcd}[row sep=tiny]
		\bullet_{R_1}^{\nu(R_1)}\arrow[dr, dash, bend right]&&  \\
		&\bullet_{R_2}^{\nu(R_2)}\arrow[r,dash,dotted]&\bullet^{\nu(R_{t-1})}_{R_{t-1}} \\
	\end{tikzcd}
	
	\qquad(b)\begin{tikzcd}[row sep=tiny]
		\bullet_{R_1}^{\nu(R_1)}\arrow[dr, dash, bend right]&&&  \\
		&\bullet_{R_2}^{\nu(R_2)}\arrow[r,dash,dotted]&\bullet^{\nu(R_{t-1})}_{R_{t-1}}\arrow[r,dash,dotted]&\bullet^{1}_{R_t}\arrow[dd,dash] \\
		&&&\\
		&&&\bullet_W^1\\
	\end{tikzcd}
	
	\qquad(c)\begin{tikzcd}[row sep=tiny]
		\bullet_{R_1}^{\nu(R_1)}\arrow[dr, dash, bend right]&&&&&  \\
		&\bullet_{R_2}^{\nu(R_2)}\arrow[r,dash,dotted]&\bullet^{d-1}_{R_t}\arrow[dr,dash,bend left]&\\
		&&&\bullet_{U}^2\arrow[r,dash]&\bullet_{W_1}^1\arrow[r,dash,dotted]&\bullet_{W_{d-3}}^1\\
	\end{tikzcd}
	
	\caption{The Enriques diagram $(E_D,\lambda_D)$.\label{diagrE}}
\end{figure}
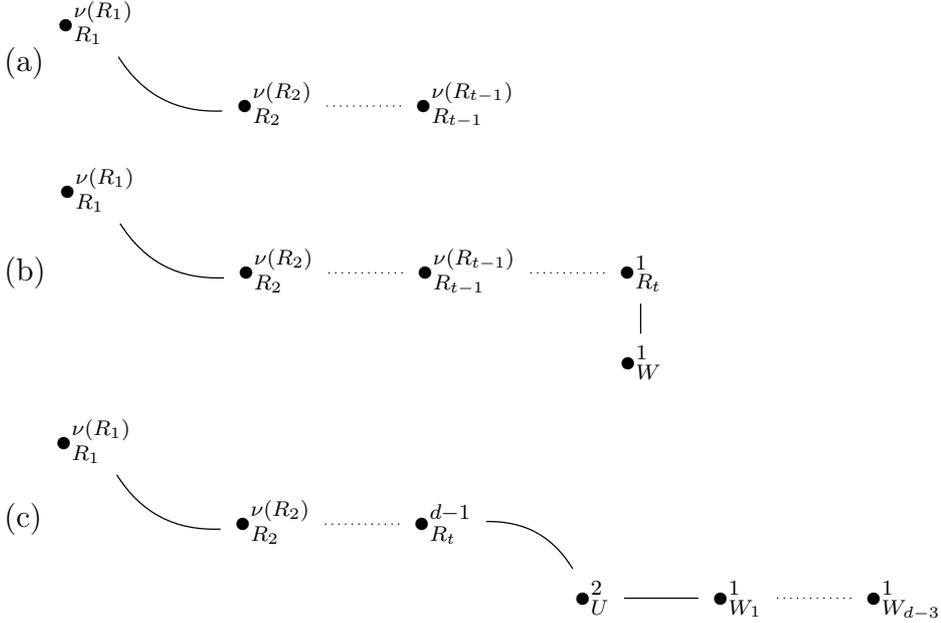

We will start from the easier part i.e. we will find the Enriques diagrams which realize these values. This theorem will be proved even for any $(D,\nu)\in H^t_d$ (not only for $(D,\nu)\in Q^t_d$).
\begin{twr}\label{latw}
	Let $d,t\in\mathbb{N}$, $dt>1$ and $(D,\nu)$ be an Enriques diagram from $H_d^t$. There exists a minimal Enriques diagram $(E_D,\lambda_D)\notin[(D,\nu)]$ such that $[(D,\nu)]$ is linear adjacent to $\left[(E_D,\lambda_D)\right]$ and
	\begin{equation}\label{wzor}\mu\left((E_D,\lambda_D)\right)=\left\{\begin{array}{ll}
			\mu((D,\nu))-1,&\mathrm{if}\hspace{1ex} d=1\\
			\mu((D,\nu))-1,&\mathrm{if}\hspace{1ex} d=2,w_D=0\\
			\mu((D,\nu))-w_D,&\mathrm{if}\hspace{1ex}d=2,w_D\neq 0\\
			\mu((D,\nu))-(d-2+w_D),&\mathrm{if}\hspace{1ex}d\geq 3\end{array}\right..\end{equation}
\end{twr}
\textit{Proof.}
The minimal diagram $(D,\nu)$ is shown in Figure \ref{diagrD}. 
We will define the diagram $(E_D,\lambda_D)$ by a modification of $(D,\nu)$. If $d=1$ we remove only the last vertex from $(D,\nu)$ (Figure \ref{diagrE}(a)) and this will be $(E_D,\lambda_D)$. If $d=2$ and $R_t$ is the root, then $E_D$ consists of only one vertex with weight $1$. If $d=2$ and $R_t$ is not the root we change the weight of the last vertex to $1$ and add one additional vertex $W$ with weight 1, so that $W\to R_t,R_{t-1}$ (Figure \ref{diagrE}(b)) and this is $(E_D,\lambda_D)$.  If $d\geq 3$ we change the weight of the last vertex to $d-1$ and add new vertices $U,W_1,\ldots,W_{d-3}$ (if $d=3$ there is no $W_i$ vertices), all proximate to $R_t$. The weights of new vertices are: $\lambda_D(U)=2$, $\lambda_D(W_i)=1$ (for $i=1,\ldots,d-3$). The proximity relation between new vertices is (Figure \ref{diagrE}(c))  
\begin{align*}
	&W_{d-3}\to W_{d-4},R_t\\
	&\ldots\\
	&W_2\to W_1,R_t\\
	&W_1\to U,R_t\\
	&U\to R_t.
\end{align*}

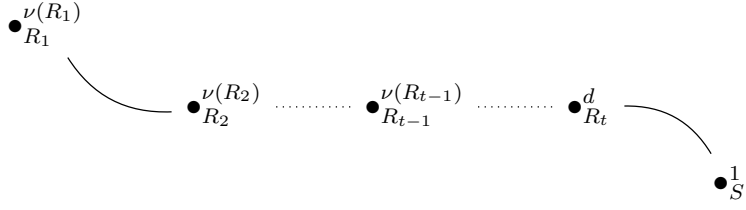
\begin{figure}[ht]
	\begin{center}
		\begin{tikzcd}[row sep=tiny]
			\bullet_{R_1}^{\nu(R_1)}\arrow[dr, dash, bend right]&&&  \\
			&\bullet_{R_2}^{\nu(R_2)}\arrow[r,dash,dotted]&\bullet^{\nu(R_{t-1})}_{R_{t-1}}\arrow[r,dash,dotted]&\bullet^d_{R_t}\arrow[dr, dash, bend left] \\
			&&&&\bullet_S^1
		\end{tikzcd}
	\end{center}
	\caption{The Enriques diagram $(D',\nu')$.\label{diagrDprim}}
\end{figure}

It is easy to check that each $(E_D,\lambda_D)$ is a minimal (and hence consistent) diagram and that $(E_D,\lambda_D)\notin[(D,\nu)]$. Moreover $(D',\nu')\geq(E_D,\lambda_D)$, where $(D',\nu')\in[(D,\nu)]$ has one additional free vertex $S$ (Figure \ref{diagrDprim}). Thus $[(D,\nu)]$ is linear adjacent to $\left[(E_D,\lambda_D)\right]$.
Now we may compute the Milnor number of $(E_D,\lambda_D)$. It is easy to notice that
$$r_{E_D}=\left\{\begin{array}{ll}
	r_D+1,&\mathrm{if}\hspace{1ex} d=1\\
	r_D-1,&\mathrm{if}\hspace{1ex} d=2,w_D=0\\
	r_D-2+w_D,&\mathrm{if}\hspace{1ex} d=2,w_D\neq 0\\ 
	r_D-d+2+w_D,&\mathrm{if}\hspace{1ex} d\geq 3\\
\end{array}\right.$$
and then after simply calculation we get (\ref{wzor}).\textrm{\hfill}$\square$

To show that the diagram from Theorem \ref{latw} realizes the maximum in (\ref{set}), it is enough to prove that for every $(D,\nu)\in Q^t_d$ all diagrams $(\tilde{D},\tilde{\nu})$ such that $[(D,\nu)]$ is linear adjacent to $[(\tilde{D},\tilde{\nu})]$ have not greater Milnor numbers than the diagram $(E_D,\lambda_D)$ constructed for $(D,\nu)$ in Theorem \ref{latw}. Of course, we may consider only $(\tilde{D},\tilde{\nu})$ which have the type different from $(D,\nu)$.  We do this in a series of lemmas in which we consecutively assume:
\begin{enumerate}
	\item\label{c1} Case - there is no subdiagram of $\tilde{D}$ isomorphic (as rooted tree with preserving shapes of edges but not weights) to $D$ (Lemma \ref{l1});
	\item\label{c2} Case - there is a subdiagram of $\tilde{D}$ isomorphic to $D$,
	\begin{enumerate}
		\item\label{sc1} Subcase - $d>2$,
		\begin{enumerate}
			\item The inequality $$\sum\limits_{P\textrm{ successor of }i^{-1}(R_t)}\min(2,\tilde{\nu}(P))+z\leq \nu(R_t)-1,$$
			where $z$ is the number of vertices proximate to $i^{-1}(R_t)$ that are not its successors, holds (Lemma \ref{l2});			
			\item The opposite inequality $$\sum\limits_{P\textrm{ successor of }i^{-1}(R_t)}\min(2,\tilde{\nu}(P))+z>\nu(R_t)-1,$$
			where $z$ is the number of vertices proximate to $i^{-1}(R_t)$ that are not its successors, holds (Lemma \ref{l4});
		\end{enumerate}
		\item\label{sc2} Subcase - $d=2$ (Lemma \ref{l5});
		\item\label{sc3} Subcase - $d=1$ (Lemma \ref{l6}).
	\end{enumerate}
\end{enumerate}	
We start with the case (\ref{c1}) that there is no subdiagram of $\tilde{D}$ isomorphic to D.
\begin{lem}\label{l1}
	Let $d,t\in\mathbb{N}$, $(D,\nu)\in Q_d^t$ and let $(\tilde{D},\tilde{\nu})$ be an arbitrary Enriques diagram such $[(D,\nu)]$ is linear adjacent to $[(\tilde{D},\tilde{\nu})]$. If there is no subdiagram of $\tilde{D}$ isomorphic to $D$, then 
			$$\mu(\tilde{D},\tilde{\nu}))\leq\left\{\begin{array}{ll}
		\mu((D,\nu))-1,&\mathrm{if}\hspace{1ex} d=1\\
		\mu((D,\nu))-1,&\mathrm{if}\hspace{1ex} d=2,w_D=0\\
		\mu((D,\nu))-w_D,&\mathrm{if}\hspace{1ex}d=2,w_D\neq 0\\
		\mu((D,\nu))-(d-2+w_D),&\mathrm{if}\hspace{1ex}d\geq 3\end{array}\right..$$
\end{lem}
\textit{Proof.} Firstly, assume that $(\tilde{D},\tilde{\nu})$ is a minimal Enriques diagram. Now, we will construct another diagram $(E,\lambda)$ such that $[(E,\lambda)]$ is linear adjacent to $[(\tilde{D},\tilde{\nu})]$ and $\mu\left((E,\lambda)\right)=\mu\left((E_D,\lambda_D)\right)$. Since $[(D,\nu)]$ is linear adjacent to $[(\tilde{D},\tilde{\nu})]$ there exist a consistent $(D',\nu')\in[(D,\nu)]$ such that $(D',\nu')\geq(\tilde{D},\tilde{\nu})$, two subdiagrams $\tilde{D}_0\subset \tilde{D}$, $D'_0\subset D'$ and an isomorphism $i:D'_0\to \tilde{D}_0$. Let $(E_D,\lambda_D)$ be the diagram from Theorem \ref{latw} constructed for $(D,\nu)$ (of course $[(D',\nu')]$ is also linear adjacent to $[(E_D,\lambda_D)]$). Since there is no subdiagram of $\tilde{D}$ isomorphic to $D$, we have $R_t\notin D_0'$ and consequently for every $P\in\tilde{D}_0$ it holds $\kappa_{\lambda}(P)=\kappa_{\nu'}(P)$ (diagrams $(D,\nu)$ and $(E,\lambda)$ are different "after $R_t$").
Then a modification of $E_D$ (analogous to the construction of $D'$ from $D$) should be made to get a diagram $(E,\lambda)\in[(E_D,\lambda_D)]$. This implies that $[(E,\lambda)]$ is linear adjacent to $[(\tilde{D},\tilde{\nu})]$, so for every singularity $f_0$ whose Enriques diagram belong to $[(E,\lambda)]$, there exists a linear deformation $(f_s)$ of $f_0$ such that the Enriques diagram of a generic element $f_s$ belongs to $[(\tilde{D},\tilde{\nu})]$ (Theorem \ref{twrmac}). Because the Milnor number is upper semi-continuous (\cite{GLS} Theorem 2.6) then for sufficiently small $s$, we have $\mu(f_s)\leq \mu(f_0)$. Therefore $\mu((\tilde{D},\tilde{\nu}))=\mu(f_s)\leq\mu(f_0)=\mu\left((E,\lambda)\right)=\mu\left((E_D,\lambda_D)\right)$.  \hfill$\square$

In the next lemmas we will consider the case (\ref{c2}) that there exists subdiagram of $\tilde{D}$ isomorphic to $D$. First, the two lemmas for the subcase (\ref{sc1}) $d>2$.

\begin{lem}\label{l2}
	Let $d,t\in\mathbb{N}$, $d\geq 3$, $(D,\nu)\in Q_t^d$ and let $(\tilde{D},\tilde{\nu})\notin[(D,\nu)]$ be an arbitrary Enriques diagram such $[(D,\nu)]$ is linear adjacent to $[(\tilde{D},\tilde{\nu})]$. If
	\begin{enumerate}
		\item there exist a subdiagram $\tilde{D}_0\subset \tilde{D}$ and an isomorphism $i:\tilde{D}_0\to D$ (not necessarily preserving the weights),
		\item $$\sum\limits_{P\textrm{ successor of }i^{-1}(R_t)}\min(2,\tilde{\nu}(P))+z\leq\nu(R_t)-1,$$
		where $z$ is the number of vertices proximate to $i^{-1}(R_t)$ that are not its successors,
	\end{enumerate} 
	then $$\mu((\tilde{D},\tilde{\nu}))\leq
	\mu((D,\nu))-(d-2+w_D).$$
\end{lem}
\textit{Proof.} We may assume that $(\tilde{D},\tilde{\nu})$ is a minimal diagram. Notice that \begin{equation}\label{nie} ord_{\tilde{\nu}}(i^{-1}(R_t))<ord_{\nu}(i^{-1}(R_t)).
\end{equation}
In fact, we prove this by induction with respect to the number of satellite vertices in $D$. Let us pass to the construction of $(E,\lambda)$ such that $[(E,\lambda)]$ is linear adjacent to $[(\tilde{D},\tilde{\nu})]$ and $\mu\left((E,\lambda)\right)\leq
\mu((D,\nu))-(d-2+w_D)$. We do this in two steps, first we construct $(E',\lambda')$ and then after some simple modification of $(E',\lambda')$ we get $(E,\lambda)$.

 Let $\{S_1,\ldots,S_m\}$ be the set of vertices proximate to $i^{-1}(R_t)$. We will construct $(E',\lambda')$. 
\begin{itemize}
	\item $E'=\{Q_1,\ldots,Q_t,U_1,\ldots,U_m\}$,
	\item $\lambda'(Q_i)=\nu(R_i)$ for $i=1,\ldots,t-1$,
	\item $\lambda'(Q_t)=\nu(R_t)-1$,
	\item $\lambda'(U_i)=\min(2,\tilde{\nu}(S_i))$ for $S_i$ that are free ($i\in\{1,\ldots,m\}$),
	\item $\lambda'(U_i)=1$ for $S_i$ that are not free ($i\in\{1,\ldots,m\}$),
	\item $Q_i\xrightarrow{E'}Q_j\Leftrightarrow R_i\xrightarrow{D}R_j$ for $i,j\in\{1,\ldots,t\}$,
	\item $U_i\xrightarrow{E'}U_j\Leftrightarrow S_i\xrightarrow{\tilde{D}}S_j$ for $i,j\in\{1,\ldots,m\}$,
	\item $S_i\xrightarrow{\tilde{D}}i^{-1}(R_k)\Rightarrow U_i\xrightarrow{E'}Q_t$ for $i\in\{1,\ldots,m\}$,
	\item $U_i\xrightarrow{E'} Q_t$ for $i=1,\ldots,m$.	
\end{itemize} 
The diagram $(E',\lambda')$ is consistent due to the second condition in the assumption. Its Milnor number can be easily estimated by $$\mu((E,\lambda))=\mu((D,\nu))-(d-2+w_D)-d^2+3d-2-x\leq \mu((D,\nu))-(d-2+w_D),$$ where $x$ is the number of successors of $i^{-1}(R_t)$ in $\tilde{D}$ with weight $1$. Because $[(D,\nu)]$ is linear adjacent to $[(\tilde{D},\tilde{\nu})]$, there exists $(D',\nu')\in [(D,\nu)]$ such that $(D',\nu')\geq (\tilde{D},\tilde{\nu})$. We can modify $(E',\lambda')$ to get $(E,\lambda)\in[(E',\lambda')]$ (analogous to the construction of $D'$ from $D$). Then  for every $P\in \tilde{D}\setminus\{i^{-1}(R_t)\}$ we have $\textrm{ord}_{\nu'}(P)\leq \textrm{ord}_{\lambda'}(P)$ and $\textrm{ord}_{\nu'}(i^{-1}(R_t))-1= \textrm{ord}_{\lambda'}i^{-1}(R_t)$. From these facts and (\ref{nie}) we get that $(E',\lambda')\geq (\tilde{D},\tilde{\nu})$. This gives that $[(E,\lambda)]$ is linear adjacent to $[(\tilde{D},\tilde{\nu})]$, so for every singularity $f_0$ whose Enriques diagram belong to $[(E,\lambda)]$, there exists a linear deformation $(f_s)$ of $f_0$ such that the Enriques diagram of a generic element $f_s$ belongs to $[(\tilde{D},\tilde{\nu})]$ (Theorem \ref{twrmac}). Because the Milnor number is upper semi-continuous (\cite{GLS} Theorem 2.6) then for sufficiently small $s$, we have $\mu(f_s)\leq \mu(f_0)$. Therefore $\mu((\tilde{D},\tilde{\nu}))\leq\mu\left((E,\lambda)\right)=\mu\left((E',\lambda')\right)\leq \mu((D,\nu))-(d-2+w_D)$. \hfill$\square$

Now, we will consider the opposite situation to the second condition in Lemma \ref{l2}.

\begin{lem}\label{l4}
	Let $d,t\in\mathbb{N}$, $d\geq 2$, $(D,\nu)\in Q_d^t$ and let $(\tilde{D},\tilde{\nu})\notin[(D,\nu)]$ be an arbitrary Enriques diagram such that $[(D,\nu)]$ is linear adjacent to $[(\tilde{D},\tilde{\nu})]$. Let us assume there exist a subdiagram $\tilde{D}_0\subset \tilde{D}$ and an isomorphism $i:\tilde{D}_0\to D$ such that $$\sum\limits_{P\textrm{ successor of }i^{-1}(R_t)}\min(2,\tilde{\nu}(P))+z>\nu(R_t)-1,$$
	where $z$ is number of vertices proximate to $i^{-1}(R_t)$ that are not its successors. Then  $\mu((\tilde{D},\tilde{\nu}))\leq
	\mu((D,\nu))-(d-2+w_D)$.
\end{lem}
\textit{Proof.} Since $\tilde{\nu}(i^{-1}(R_t))\leq \nu(R_t)$ and  $$\sum\limits_{P\textrm{ successor of }i^{-1}(R_t)}\min(2,\tilde{\nu}(P))+z\leq \tilde{\nu}(i^{-1}(R_t)),$$ we get $$\nu(R_t)\leq\sum\limits_{P\textrm{ successor of }i^{-1}(R_t)}\min(2,\tilde{\nu}(P))+z\leq\tilde{\nu}(i^{-1}(R_t))\leq\nu(R_t).$$ Then we get hence the equality $\nu(R_t)=\tilde{\nu}(i^{-1}(R_t))$. Because $(\tilde{D},\tilde{\nu})\notin[(D,\nu)]$, in $(\tilde{D},\tilde{\nu})$ all the successors of $i^{-1}(R_t)$ have weight $2$ at most. If after them there is a vertex of weight $2$, it has to be free. So after the $i^{-1}(R_t)$ we can have a "new branch" with vertices of weight $2$. The length of such "new branch" is limited by $\textrm{ord}_{\nu}(R_t)-\textrm{ord}_{\tilde{\nu}}(i^{-1}(R_t))$.

Moreover, since $\nu(R_t)=\tilde{\nu}(i^{-1}(R_t))$ also $\nu(R_j)=\tilde{\nu}(i^{-1}(R_j))$ for $j=k_0+1,\ldots,t$, where $R_{t_0}$ is the last free vertex in $(D,\nu)$. 

The Milnor number of $(\tilde{D},\tilde{\nu})$ can be estimated by:
\begin{multline*}
	\mu((\tilde{D},\tilde{\nu}))\leq \mu((D,\nu))+\\ \sum\limits_{j=1}^{k_0}\left(
	\tilde{\nu}(i^{-1}(R_j))\left(\tilde{\nu}(i^{-1}(R_j))-1\right)-\nu(R_i)\left(\nu(R_i)-1\right)
	\right)+\\d\left(\textrm{ord}_{\nu}(R_t)-\textrm{ord}_{\tilde{\nu}}(i^{-1}(R_t))\right)+(\nu(R_1)-\tilde{\nu}(i^{-1}(R_1)))\leq\\ \mu((D,\nu))-(d-2+w_D).
\end{multline*}
\hfill$\square$

In the next lemma we consider the subcase (\ref{sc2}) $d=2$.

\begin{lem}\label{l5}
	Let $k\in\mathbb{N}$, $(D,\nu)\in Q_2^t$ and let $(\tilde{D},\tilde{\nu})\notin[(D,\nu)]$ be an arbitrary minimal Enriques diagram such that $[(D,\nu)]$ is linear adjacent to $[(\tilde{D},\tilde{\nu})]$. Then  $$\mu((\tilde{D},\tilde{\nu}))\leq\left\{\begin{array}{ll}
		\mu((D,\nu))-1,&\mathrm{if}\hspace{1ex}w_D=0\\
		\mu((D,\nu))-w_D,&\mathrm{if}\hspace{1ex}w_D\neq 0\end{array}\right.$$
\end{lem}
\textit{Proof.} 
If $w_D=0$ then the only diagram $(\tilde{D},\tilde{\nu})\notin[(D,\nu)]$ such that $[(D,\nu)]$ is linear adjacent to $[(\tilde{D},\tilde{\nu})]$ is $(E_D,\lambda_D)$. Then $\mu((\tilde{D},\tilde{\nu}))=\mu((E_D,\lambda_D))=\mu((D,\nu))-1$.

Let assume that $w_D\neq 0$. If there is no subdiagram of $\tilde{D}$ isomorphic to $D$ we can apply Lemma \ref{l1}. If there exist subdiagrams $\tilde{D}_0\subset \tilde{D}$, $D_0\subset D$ and an isomorphism $i:\tilde{D}_0\to D$, then $$\sum\limits_{P\textrm{ successor of }i^{-1}(R_t)}\min(2,\tilde{\nu}(P))+z=\nu(R_t),$$
where $z$ is the number of vertices proximate to $i^{-1}(R_t)$ that are not its successors. Then from Lemma \ref{l4} we get $\mu((\tilde{D},\tilde{\nu}))\leq\mu((D,\nu))-w_D$.\hfill$\square$

The last lemma is for the last subcase (\ref{sc3}) $d=1$ and it is easy to prove.
\begin{lem}\label{l6}
	Let $k\in\mathbb{N}$, $(D,\nu)\in Q_1^t$ and let $(\tilde{D},\tilde{\nu})\notin[(D,\nu)]$ be an arbitrary minimal Enriques diagram such that $[(D,\nu)]$ is linear adjacent to $[(\tilde{D},\tilde{\nu})]$. If there exist subdiagrams $\tilde{D}_0\subset \tilde{D}$, $D_0\subset D$ and an isomorphism $i:\tilde{D}_0\to D$, then $\mu((\tilde{D},\tilde{\nu}))\leq
	\mu((D,\nu))-1$.
\end{lem}

Now, we can formulate the main result (that indeed the diagram $(E_D,\lambda_D)$ from Theorem \ref{latw} realizes the maximum in (\ref{set})). This theorem is a consequence of previous lemmas.

\begin{twr}\label{tru}	Let $d,t\in\mathbb{N}$, $dt>1$, $(D,\nu)\in Q_d^t$ and let $(\tilde{D},\tilde{\nu})\notin[(D,\nu)]$ be an arbitrary Enriques diagram such that  $[(D,\nu)]$ is linear adjacent to $[(\tilde{D},\tilde{\nu})]$. Then 	$$\mu((\tilde{D},\tilde{\nu}))\leq\left\{\begin{array}{ll}
		\mu((D,\nu))-1,&\mathrm{if}\hspace{1ex}d=1\\
		\mu((D,\nu))-1,&\mathrm{if}\hspace{1ex}d=2,w_D=0\\
		\mu((D,\nu))-w_D,&\mathrm{if}\hspace{1ex}d=2,w_D\neq 0\\
		\mu((D,\nu))-(d-2+w_D),&\mathrm{if}\hspace{1ex}d\geq 3\end{array}\right..$$
\end{twr}

\section{Formula for jump of the Milnor number of quasihomogeneous singularity for linear deformations}
In this section we apply Theorems \ref{latw} and \ref{tru} to the Enriques diagrams of a quasihomogeneous singularity.

As a consequence of these theorems and the construction of the Enriques diagrams of quasihomogeneous singularities we can formulate the following facts.
\begin{twr}\label{quasitwr}
	For any quasihomogeneous singularity 
	$f_0$ of form (\ref{quasi}) the jump of Milnor number of $f_0$ for linear deformations is
	\begin{equation}\label{main}\lambda^{lin}(f_0)=\left\{\begin{array}{ll}
			1,&\mathrm{if}\hspace{1ex}d=1\\
			1,&\mathrm{if}\hspace{1ex}d=2,w_{f_0}=0\\
			w_{f_0},&\mathrm{if}\hspace{1ex}d=2,w_{f_0}\neq 0\\
			d-2+w_{f_0},&\mathrm{if}\hspace{1ex}d\geq 3\end{array}\right.,\end{equation}
	where:\begin{itemize}
		\item if $p=q$ then $d=k+l+p$ and $w_{f_0}=0$,
		\item if $p\neq q$ and $p\vert q$ then $d=k+p$ and $w_{f_0}=1$,
		\item if $p\neq q$ and $p\doesnotdivide q$ then $d=\textrm{GCD}(p,q)$ and $w_{f_0}=2$.
	\end{itemize}
\end{twr}
\textit{Proof.} Let $f_0$ be a quasihomogeneous singularity and $(D,\nu)$ its Enriques diagram. From Theorem \ref{latw} there exists diagram $(E_D,\lambda_D)\notin[(D,\nu)]$ such that $[(D,\nu)]$ is linear adjacent to $\left[(E_D,\lambda_D)\right]$ and $$\mu\left((E_D,\lambda_D)\right)=\left\{\begin{array}{ll}
	\mu((D,\nu))-1,&\mathrm{if}\hspace{1ex} d=1\\
	\mu((D,\nu))-1,&\mathrm{if}\hspace{1ex} d=2,w_D=0\\
	\mu((D,\nu))-w_D,&\mathrm{if}\hspace{1ex}d=2,w_D\neq 0\\
	\mu((D,\nu))-(d-2+w_D),&\mathrm{if}\hspace{1ex}d\geq 3\end{array}\right..$$ 

Since $(E_D,\lambda_D)$ is minimal Theorem \ref{bije} and Theorem \ref{twrmac} give 
\begin{equation}\label{ine}\lambda^{lin}(f_0)\leq\left\{\begin{array}{ll}
		1,&\mathrm{if}\hspace{1ex}d=1\\
		1,&\mathrm{if}\hspace{1ex}d=2,w_D=0\\
		w_D,&\mathrm{if}\hspace{1ex}d=2,w_D\neq 0\\
		d-2+w_D,&\mathrm{if}\hspace{1ex}d\geq 3\end{array}\right..\end{equation}

From Theorem \ref{tru} for any Enriques diagram $(\tilde{D},\tilde{\nu})\notin[(D,\nu)]$  such that $[(D,\nu)]$ is linear adjacent to $[(\tilde{D},\tilde{\nu})]$ we have $\mu((\tilde{D},\tilde{\nu}))\leq\mu\left((E_D,\lambda_D)\right)$. It gives the opposite inequality in (\ref{ine}) and as a consequence we get (\ref{main}), because $w_D=w_{f_0}$.
\hfill$\square$

Taking into account the simple characterization (\ref{w}) of $w_D$, we get a more effective formula.
\begin{wn}\label{wn1}
	Let $f_0$ be a quasihomogeneous singularity  of form (\ref{quasi}). Then
	\begin{enumerate}
		\item If $p=q$ i.e. $f_0$ is a homogeneous singularity then
		$$\lambda^{lin}(f_0)=\left\{\begin{array}{ll}
			1,&\mathrm{if}\hspace{1ex}k+l+p=2\\
			k+l+p-2,&\mathrm{if}\hspace{1ex}k+l+p\geq 3\end{array}\right..$$
		\item If $p\neq q$ and $p\vert q$ then
		$$\lambda^{lin}(f_0)=\left\{\begin{array}{ll}
			1,&\mathrm{if}\hspace{1ex}p+k\leq 2\\
			p+k-1,&\mathrm{if}\hspace{1ex}p+k\geq 3\end{array}\right..$$		
		\item If $p\neq q$ and $p\doesnotdivide q$ then
		$$\lambda^{lin}(f_0)=GCD(p,q).$$		
	\end{enumerate}
\end{wn}

If we consider only the "standard" quasihomogeneous singularities i.e. $k=l=0$ in (\ref{quasi}), we get a very simple formula for the jump.
\begin{wn}
	Let $f_0$ be a quasihomogeneous singularity defined in (\ref{quasi}) and $k=l=0$. Then
	\begin{enumerate}
		\item If $p=q$ i.e. $f_0$ is a homogeneous singularity then
		$$\lambda^{lin}(f_0)=\left\{\begin{array}{ll}
			1,&\mathrm{if}\hspace{1ex}p=2\\
			p-2,&\mathrm{if}\hspace{1ex}p\geq 3\end{array}\right..$$
		\item If $p\neq q$ then 
		$$\lambda^{lin}(f_0)=\left\{\begin{array}{ll}
		p-1,&\mathrm{if}\hspace{1ex}p\vert q\\
		GCD(p,q),&\mathrm{if}\hspace{1ex}p\doesnotdivide q\end{array}\right..$$				
	\end{enumerate}
\end{wn}
\begin{przy}\label{prz3}
	Let's consider the singularity from Example \ref{prz2} i.e. $f_0(x,y)=x^6+y^9$. Its minimal Enriques diagram is shown in Figure \ref{rprzy}(a). The minimal Enriques diagram realizing the $\lambda^{lin}(f_0)$ (constructed in Theorem \ref{latw}) is shown in Figure \ref{rprzy}(b). A linear deformation having this diagram is $f_s(x,y)=f_0(x,y)+sxy(y^3+x^2)^2$. 
	\begin{figure}[ht]
		\begin{center}
			(a)\begin{tikzcd}[row sep=tiny]
			&\bullet^3\arrow[dd,dash]\\
			&\\
			\bullet^6\arrow[uur,dash,bend left]&\bullet^3\\
			\end{tikzcd}
			(b)\begin{tikzcd}[row sep=tiny]
				&\bullet^3\arrow[dd,dash]&\\
				&&\\
				\bullet^6\arrow[uur,dash,bend left]&\bullet^2\arrow[dr,dash,bend right]&\\
				&&\bullet^2\\
			\end{tikzcd}		
		\end{center}
		\caption{The minimal Enriques diagrams of $f_0$ and $f_s$\label{rprzy}}
	\end{figure}
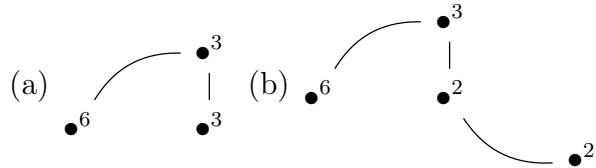
\end{przy}
\begin{uwa}
	It is not an easy task to write down an explicit formula of the deformation from the constructed Enriques diagram. Obviously, in specific case it can be done (as in Example \ref{prz3}).
\end{uwa}
As a corollary we give a formula for the jump of Milnor number for \textbf{semi-quasiho\-mogeneous singularities} i.e. singularities of the form $f_0=f_0'+g,$ where $f_0'$ is a quasihomogeneous singularity with respect to some weights $(w_x,w_y)$ and $\textrm{ord}_{(w_x,w_y)} g>\textrm{ord}_{(w_x,w_y)} f_0'$. 
\begin{wn}\label{semi}
	For any semi-quasihomogeneous singularity $f_0$
	$$\lambda^{lin}(f_0)=\lambda^{lin}(f_0').$$
\end{wn}
\textit{Proof.} It suffices to notice that Enriques diagrams of $f_0$ and $f_0'$ have the same type. \hfill$\square$

\textbf{Acknowledgments}

I would like to thank Tadeusz Krasi\'nski for his support in writing this article and Szymon Brzostowski for his editorial notes.
\bibliographystyle{alpha}
\bibliography{biblio}

\begin{flushright}\footnotesize
\textsc{Aleksandra Zakrzewska\\Faculty of Mathematics and Computer Science\\ University of Lodz\\ul. Banacha 22\\90-238 Lodz, Poland\\aleksandra.zakrzewska@wmii.uni.lodz.pl}
\end{flushright}

\end{document}